\documentstyle[11pt]{article}
\normalbaselines
\begin{document}
\bibliographystyle{unsrt}

\def\bea*{\begin{eqnarray*}}
\def\eea*{\end{eqnarray*}}
\def\ba{\begin{array}}
\def\ea{\end{array}}
\count1=1
\def\be{\ifnum \count1=0 $$ \else \begin{equation}\fi}
\def\ee{\ifnum\count1=0 $$ \else \end{equation}\fi}
\def\ele(#1){\ifnum\count1=0 \eqno({\bf #1}) $$ \else \label{#1}\end{equation}\fi}
\def\req(#1){\ifnum\count1=0 {\bf #1}\else \ref{#1}\fi}
\def\bea(#1){\ifnum \count1=0   $$ \begin{array}{#1}
\else \begin{equation} \begin{array}{#1} \fi}
\def\eea{\ifnum \count1=0 \end{array} $$
\else  \end{array}\end{equation}\fi}
\def\elea(#1){\ifnum \count1=0 \end{array}\label{#1}\eqno({\bf #1}) $$
\else\end{array}\label{#1}\end{equation}\fi}
\def\cit(#1){
\ifnum\count1=0 {\bf #1} \cite{#1} \else 
\cite{#1}\fi}
\def\bibit(#1){\ifnum\count1=0 \bibitem{#1} [#1    ] \else \bibitem{#1}\fi}
\def\ds{\displaystyle}
\def\hb{\hfill\break}
\def\comment#1{\hb {***** {\em #1} *****}\hb }

\newcommand{\TZ}{\hbox{\bf T}}
\newcommand{\MZ}{\hbox{\bf M}}
\newcommand{\ZZ}{\hbox{\bf Z}}
\newcommand{\NZ}{\hbox{\bf N}}
\newcommand{\RZ}{\hbox{\bf R}}
\newcommand{\CZ}{\,\hbox{\bf C}}
\newcommand{\PZ}{\hbox{\bf P}}
\newcommand{\QZ}{\hbox{\bf Q}}
\newcommand{\HZ}{\hbox{\bf H}}
\newcommand{\EZ}{\hbox{\bf E}}
\newcommand{\GZ}{\,\hbox{\bf G}}
\newcommand{\DZ}{\, \hbox{\bf D}}

\newtheorem{theorem}{Theorem}
\newtheorem{lemma}{Lemma}
\newtheorem{proposition}{Proposition}
\newtheorem{corollary}{Corollary}

\vbox{\vspace{38mm}}
\begin{center}
{\LARGE \bf On Branching Indices of Affine A-D-E Diagrams : A Geometrical Characterization by Kleinian Singularities}\\[10 mm]  
Shi-shyr Roan 
\\[2mm]
{\it Institute of Mathematics, Academia Sinica \\ 
Taipei , Taiwan \\ (e-mail: maroan@ccvax.sinica.edu.tw)} 
\\[35mm]
\end{center}

\begin{abstract}
The exceptional configuration of the minimal resolution $\widehat{S}_G $ of a Kleinian quotient surface $S_G (:= \CZ^2/G)$ is depicted by a $A$-$D$-$E$ Coxeter-Dynkin diagram. In this article, we show that branching indices of the affine $A$-$D$-$E$ diagram is  geometrically characterized by a certain special function $F$ of $S_G$ as the multiplicities of its divisor components in $\widehat{S}_G$, a version parallel to the elliptic fibration near certain types of simple singular fibers in Kodaira's elliptic surface theory. We further obtain the uniqueness property of the function $F$ (modular local units) among all local functions in $S_G$ near the singular point whose divisors in  $\widehat{S}_G $ display the affine $A$-$D$-$E$ diagram configuration.
\end{abstract}

\par \vspace{5mm} \noindent
\vfill
\eject

\section{Introduction}
In his study of the general theory of elliptic surfaces, K. Kodaira obtained the complete structure of singular fibers over a non-singular curve \cite{Kod}. 
In the list of all singular fiber structures, there are eight simple types, labeled by ${\rm I}_b, {\rm I}_b^*$, ${\rm II}, {\rm II}^*$, ${\rm III}, {\rm III}^*$, ${\rm IV}, {\rm IV}^*$. Among which, the types ${\rm I}_b \ (b \geq 2)$, ${\rm I}_b^*, {\rm II}^*, {\rm III}^*$, ${\rm IV}^*$, are special in that the irreducible components of the singular fiber are all $(-2)$-$\PZ^1$ curves with normal crossing. Furthermore, the fibration structure of the smooth surface near a such singular fiber can be represented by a marked graph, consisting of lines attached with some positive integers, to describe the singular fiber-divisor in the surface. The  structure can equally be represented by its dual graph with a branching index attached to each node. The resulting dual graphs are the affine $A$-$D$-$E$ Coxeter-Dynkin diagrams, described by the following correspondence\footnote{In this correspondence, the irreducible component $\Theta_0$ of the singular fiber in \cite{Kod} corresponds to $-\theta$ of the affine Coxeter-Dynkin diagram.}:
\bea(clc)
{\rm Kodaira \ singular \ elliptic \ fiber \ type } & \stackrel{\rm dual}{\Longleftrightarrow} & {\rm Affine \ Coxeter-Dynkin \ diagram}  \\[2mm] 
{\rm I}_{r+1} & \longleftrightarrow &  \widetilde{A}_r \ \ \ (r \geq 1) \ , \\
{\rm I}^*_b \ \ & \longleftrightarrow & \ \widetilde{D}_{4+b} \ (b \geq 0) \ ,  \\
{\rm IV}^* & \longleftrightarrow &   \widetilde{E}_6 \ \ \ \ , \ \ \ \ \ \ \ \ \ \\
{\rm III}^* & \longleftrightarrow &  \widetilde{E}_7 \ \ \ \ , \ \ \ \ \ \ \ \ \ \\
{\rm II}^* & \longleftrightarrow & \widetilde{E}_8 \ \ \ \ . \ \ \ \ \ \ \ \ \ 
\elea(KoADE)
where the affine $A$-$D$-$E$ diagrams are depicted graphically as follows: 
$$
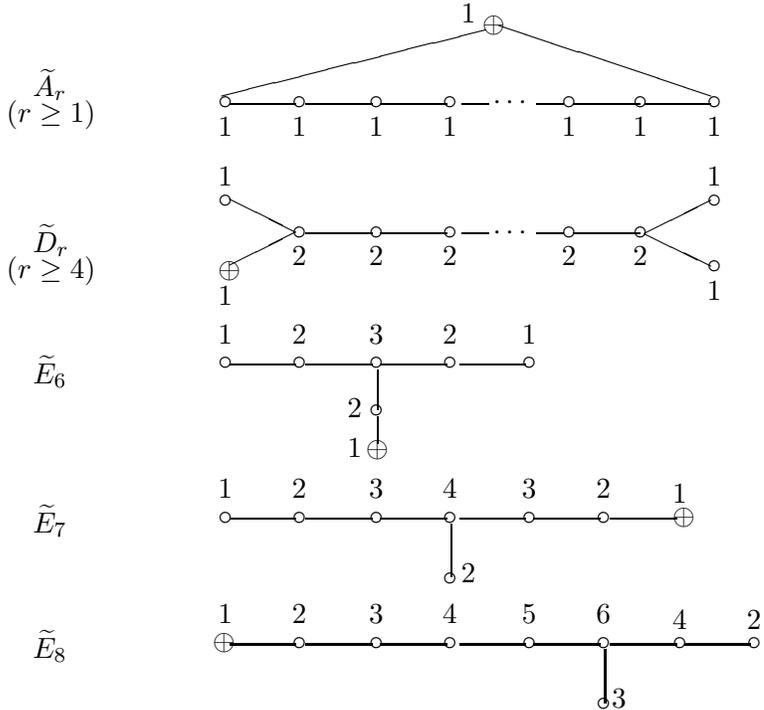
\put(-200, 5) { \shortstack{$\widetilde{A}_r$}}\\
\put(-210, -5) { \shortstack{$(r \geq 1)$}}\\
\put(-125, 5){\line(4,1){100}} \\
\put(-130, 0) { \shortstack{$\circ$}} \\
\put(-130, -10) { \shortstack{$1$}} \\
\put(-122, 2){\line(1,0){25}} \\
\put(-102, 0) { \shortstack{$\circ$}} \\
\put(-102, -10) { \shortstack{$1$}} \\
\put(-93, 2){\line(1,0){25}} \\
\put(-73, 0) { \shortstack{$\circ$}} \\
\put(-73, -10) { \shortstack{$1$}} \\
\put(-65, 2){\line(1,0){25}} \\
\put(-45, 0) { \shortstack{$\circ$}} \\
\put(-45, -10) { \shortstack{$1$}} \\
\put(-34, 2){\line(1,0){10}} \\
\put(-26, 0) { \shortstack{$\cdots$}} \\
\put(-6, 2){\line(1,0){10}} \\
\put(0, 0) { \shortstack{$\circ$}} \\
\put(0, -10) { \shortstack{$1$}} \\
\put(8, 2){\line(1,0){25}} \\
\put(27, 0) { \shortstack{$\circ$}} \\
\put(27, -10) { \shortstack{$1$}} \\
\put(35, 2){\line(1,0){25}} \\
\put(55, 0) { \shortstack{$\circ$}} \\
\put(55, -10) { \shortstack{$1$}} \\
\put(-30, 29) { \shortstack{$\oplus$}} \\
\put(-38, 33) { \shortstack{$1$}} \\
\put(60, 5){\line(-3,1){80}} \\
$$
$$
\put(-200, -5) { \shortstack{$\widetilde{D}_r$}}\\
\put(-210, -15) { \shortstack{$(r \geq 4)$}}\\
\put(-130, -15) { \shortstack{$\oplus$}} \\
\put(-130, -25) { \shortstack{$1$}} \\
\put(-122, -10){\line(2,1){25}} \\
\put(-130, 12) { \shortstack{$\circ$}} \\
\put(-130, 20) { \shortstack{$1$}} \\
\put(-122, 15){\line(2,-1){25}} \\
\put(-102, 0) { \shortstack{$\circ$}} \\
\put(-102, -10) { \shortstack{$2$}} \\
\put(-93, 2){\line(1,0){25}} \\
\put(-73, 0) { \shortstack{$\circ$}} \\
\put(-73, -10) { \shortstack{$2$}} \\
\put(-65, 2){\line(1,0){25}} \\
\put(-45, 0) { \shortstack{$\circ$}} \\
\put(-45, -10) { \shortstack{$2$}} \\
\put(-34, 2){\line(1,0){10}} \\
\put(-26, 0) { \shortstack{$\cdots$}} \\
\put(-6, 2){\line(1,0){10}} \\
\put(0, 0) { \shortstack{$\circ$}} \\
\put(0, -10) { \shortstack{$2$}} \\
\put(8, 2){\line(1,0){25}} \\
\put(27, 0) { \shortstack{$\circ$}} \\
\put(27, -10) { \shortstack{$2$}} \\
\put(35, 2){\line(2,1){25}} \\
\put(55, 12) { \shortstack{$\circ$}} \\
\put(55, 20) { \shortstack{$1$}} \\
\put(35, 2){\line(2,-1){25}} \\
\put(55, -13) { \shortstack{$\circ$}} \\
\put(55, -23) { \shortstack{$1$}} \\
$$
$$
\put(-200, -5) { \shortstack{$\widetilde{E}_6$}}\\
\put(-130, 0) { \shortstack{$\circ$}} \\
\put(-130, 10) { \shortstack{$1$}} \\
\put(-122, 2){\line(1,0){25}} \\
\put(-102, 0) { \shortstack{$\circ$}} \\
\put(-102, 10) { \shortstack{$2$}} \\
\put(-93, 2){\line(1,0){25}} \\
\put(-73, 0) { \shortstack{$\circ$}} \\
\put(-73, 10) { \shortstack{$3$}} \\
\put(-65, 2){\line(1,0){25}} \\
\put(-45, 0) { \shortstack{$\circ$}} \\
\put(-45, 10) { \shortstack{$2$}} \\
\put(-35, 2){\line(1,0){25}} \\
\put(-15, 0) { \shortstack{$\circ$}} \\
\put(-15, 10) { \shortstack{$1$}} \\
\put(-66, 0){\line(0,-1){15}} \\
\put(-73, -18) { \shortstack{$\circ$}} \\
\put(-81, -18) { \shortstack{$2$}} \\
\put(-66, -18){\line(0,-1){10}} \\
\put(-74, -33) { \shortstack{$\oplus$}} \\
\put(-81, -33) { \shortstack{$1$}} \\
$$
$$
\put(-200, -3) { \shortstack{$\widetilde{E}_7$}}\\
\put(-130, 0) { \shortstack{$\circ$}} \\
\put(-130, 10) { \shortstack{$1$}} \\
\put(-122, 2){\line(1,0){25}} \\
\put(-102, 0) { \shortstack{$\circ$}} \\
\put(-102, 10) { \shortstack{$2$}} \\
\put(-93, 2){\line(1,0){25}} \\
\put(-73, 0) { \shortstack{$\circ$}} \\
\put(-73, 10) { \shortstack{$3$}} \\
\put(-65, 2){\line(1,0){25}} \\
\put(-45, 0) { \shortstack{$\circ$}} \\
\put(-45, 10) { \shortstack{$4$}} \\
\put(-35, 2){\line(1,0){25}} \\
\put(-15, 0) { \shortstack{$\circ$}} \\
\put(-15, 10) { \shortstack{$3$}} \\
\put(-38, 0){\line(0,-1){20}} \\
\put(-7, 2){\line(1,0){25}} \\
\put(13, 0) { \shortstack{$\circ$}} \\
\put(13, 10) { \shortstack{$2$}} \\
\put(22, 2){\line(1,0){25}} \\
\put(-45, -23) { \shortstack{$\circ$}} \\
\put(-38, -22) { \shortstack{$2$}} \\
\put(42, 0) { \shortstack{$\oplus$}} \\
\put(42, 8) { \shortstack{$1$}} \\
$$
$$
\put(-200, -3) { \shortstack{$\widetilde{E}_8$}}\\
\put(-132, 0) { \shortstack{$\oplus$}} \\
\put(-130, 10) { \shortstack{$1$}} \\
\put(-122, 2){\line(1,0){25}} \\
\put(-102, 0) { \shortstack{$\circ$}} \\
\put(-102, 10) { \shortstack{$2$}} \\
\put(-93, 2){\line(1,0){25}} \\
\put(-73, 0) { \shortstack{$\circ$}} \\
\put(-73, 10) { \shortstack{$3$}} \\
\put(-65, 2){\line(1,0){25}} \\
\put(-45, 0) { \shortstack{$\circ$}} \\
\put(-45, 10) { \shortstack{$4$}} \\
\put(-35, 2){\line(1,0){25}} \\
\put(-15, 0) { \shortstack{$\circ$}} \\
\put(-15, 10) { \shortstack{$5$}} \\
\put(13, 10) { \shortstack{$6$}} \\
\put(20, 0){\line(0,-1){20}} \\
\put(13, -23) { \shortstack{$\circ$}} \\
\put(19, -22) { \shortstack{$3$}} \\
\put(-7, 2){\line(1,0){25}} \\
\put(13, 0) { \shortstack{$\circ$}} \\
\put(22, 2){\line(1,0){25}} \\
\put(42, 0) { \shortstack{$\circ$}} \\
\put(42, 8) { \shortstack{$4$}} \\
\put(50, 2){\line(1,0){25}} \\
\put(70, 0) { \shortstack{$\circ$}} \\
\put(70, 8) { \shortstack{$2$}} \\
$$
In the above diagrams, the $\circ$'s denote fundamental roots in a Coxeter-Dynkin diagram, $\oplus = -\theta$ (the negative of maximal root), and the integer close to a node indicates the corresponding branching coefficient. Note that the total sum of roots with branching coefficients is equal to zero. 

On the other hand, in his work on the invariant theory of regular solids in $\RZ^3$ around 1872, F. Klein obtained the structure of surface singularities of the quotient space $S_G \ (:= \CZ^2/G)$  for a finite non-trivial subgroup $G$ of ${\rm SL}_2 ( \CZ )$ \cite{Kl}. The orbifold $S_G$ can be realized as a hypersurface in $\CZ^3$ with the isolated singularity at the origin, denoted by $o$ throughout the paper.  
It is known that the minimal resolution of $S_G$,
\be
\pi : \ \widehat{S}_G \longrightarrow S_G \ ,
\ele(resl) 
has the trivial canonical bundle, and the exceptional configuration $\pi^{-1}(o)$ is represented by  the $A$-$D$-$E$  Coxeter-Dynkin diagram \cite{Ar, Br, HH, Li}. The classification of Kleinian singularities 
is given by the following table:

$$
\put(-220, 20){\line(0, -1){130}}
\put(200, 20){\line(0, -1){130}}
\put(-220, 20){\line(1, 0){420}}
\put(-220, 0){\line(1, 0){420}}
\put(-190, 20){\line(0, -1){130}}
\put(-90, 20){\line(0, -1){130}}
\put(-40, 20){\line(0, -1){130}}
\put(65, 20){\line(0, -1){130}}
\put(-220, 5) { \shortstack{Type}}
\put(-165, 5){ \shortstack{Group \ G }}
\put(-90, 5){ \shortstack{ $|G|$ }}
\put(-215, -20){\shortstack{$A_r$}}
\put(-165, -20){ \shortstack{ Cyclic }}
\put(-30, -20){ \shortstack{ $2, r+1, r+1$ }}
\put(80, -20){ \shortstack{ 
$X^{r+1} - YZ$ }}
\put(-215, -40){\shortstack{$D_r$}}
\put(-190, -40){ \shortstack{ Binary Dihedral }}
\put(-32, -40){ \shortstack{ $4,2(r-2),2(r-1)$ }}
\put(80, -40){ \shortstack{ 
$X^{r-1} - XY^2 + Z^2$ }}
\put(80, -60){ \shortstack{ 
$X^4 + Y^3 + Z^2$ }}
\put(80, -80){ \shortstack{ 
$X^3Y + Y^3 + Z^2$ }}
\put(80, -100){ \shortstack{ 
$X^5 + Y^3 + Z^2$ }}
\put(-215, -60){\shortstack{$E_6$}}
\put(-215, -80){\shortstack{$E_7$}}
\put(-215, -100){\shortstack{$E_8$}}
\put(-190, -60){ \shortstack{ Binary Tetrahedral }}
\put(-190, -80){ \shortstack{ Binary Octahedral }}
\put(-190, -100){ \shortstack{ Binary Icosahedral }}
\put(-80, -60){ \shortstack{24 }}
\put(-80, -80){ \shortstack{ 48 }}
\put(-80, -100){ \shortstack{ 120 }}
\put(-85, -40){ \shortstack{ $4(r-2)$ }}
\put(-85, -20){ \shortstack{ $r+1$ }}
\put(-25, -60){ \shortstack{ 6, 8, 12 }}
\put(-25, -80){ \shortstack{ 8, 12, 18}}
\put(-25, -100){ \shortstack{ 12, 20, 30 }}
\put(-40, 5){ \shortstack{Degree of Invariants}}
\put(80, 5){ \shortstack{Relation of Invariants}}
\put(-220, -110){\line(1, 0){420}}
\put(-120, -130){ \shortstack{ Table 1: Classification of Kleinian Singularities }}
$$
where $X, Y, Z$ are $G$-invariant homogenous polynomials in the coordinates of $\CZ^2$ with the degrees indicated respectively by those in the above table. The hypersurface structure of $S_G$ in $\CZ^3$ is defined as zeros of the invariants' relation. The configuration of exceptional divisors in $\pi^{-1}(o) \subset \widehat{S}_G$ can be realized by its dual graph, which is a $A$-$D$-$E$ diagram, i.e. the one obtained by the affine Coxeter-Dynkin diagram by omitting $\oplus$ and  branching indices. The $A$-$D$-$E$ geometrical nature of $\widehat{S}_G$ in connection with the representation theory of $G$ was displayed in \cite{GSV} through the McKay correspondence \cite{Mc}, in which  a remarkable observation was discovered by  J. McKay about a natural dual-pairing between conjugacy classes and irreducible representations of a Klein group in accord with its affine Coxeter-Dynkin diagram, where the branching indices stand a simple dimensionality-explanation of representation spaces.  
Nevertheless, what still remains unclear is the geometrical interpretation in the McKay correspondence about {\it branching indices}  of the affine Coxeter-Dynkin diagram.
By Kodaira's elliptic surface theory and the correspondence (\req(KoADE)), the elliptic projection should in principle provide a local analytic function of $S_G$ near the singular point $o$ such that its divisor in $\widehat{S}_G$ via the morphism (\req(resl)) produces pictorially the affine Coxeter-Dynkin diagram. By which, the existence of a local function of $S_G$ near $o$ with the affine $A$-$D$-$E$ diagram property for its divisor in $\widehat{S}_G$ has been a well-known fact since 1960's through the Kodaira's theory on elliptic surfaces\footnote{In the surface-singularity theory, there has been an extensive study with progresses made on (isolated) rational singularities, and the multiplicity 2 case is the class of Kleinian singularities in this article. The "fundamental cycle" in the minimal resolution of a rational singularity \cite{Ar6} plays the role of the affine $A$-$D$-$E$ diagrams, but without the data $\oplus$, associated to Kleinian singularities. The structure of fundamental cycles for rational singularities and related topics has been known since 1970's (see, e.g., \cite{W} and references therein).}.   
However the explicit form of such function(s) in terms of invariants in Table 1 remains to be discovered. It is the objective of this paper to show that one can indeed determine all those functions through orbifold geometry of Kleinian quotient singularities. 
We shall provide the explicit quantitative answer to questions about functions of $S_G$ with the affine $A$-$D$-$E$ diagram  property by using methods in toric geometry and group representation theory. The main results of this article are stated in the following two theorems.
\begin{theorem}\label{thm:BrI}
Let $X, Y, Z$ be the invariants in ${\rm Table \ 1}$. We define the following function $F$ of $S_G$,
$$
F := \left\{ \begin{array}{ll} X & {\rm for} \ G = A_r, E_6, E_7, E_8 \ ; \\
X + c Y \ \ ( c \in \CZ \setminus \{ 0, \pm 1 \}) & {\rm for} \ G = D_4 ; \\
X + c Y \ \ ( c \in \CZ \setminus \{ 0 \} ) & {\rm for} \ G = D_r \ ( r \geq 5 ) \ .  \end{array} \right.
$$
Then the divisor ${\rm div}(\pi^*F)$ in $\widehat{S}_G$ defined by $\pi^*F=0$ near the exceptional set is described by its corresponding affine $A$-$D$-$E$ Coxeter-Dynkin diagram (with branching indices). Furthermore, the $\oplus$ in the diagram represents the open affine curve which is the proper transform of $(F=0)$ in $S_G$. 
\end{theorem}
Among the $A$-type groups, $A_1$ is special in that the basic invariants $X, Y, Z$ are all of degree two. One can achieve the same expression of invariants' relation for $A_1$ in Table 1 by some suitable linear changes of variables in $X, Y, Z$, then the curve corresponding to $\oplus$ varies accordingly. Thus, for the problem  of determining all local functions in $S_G$ near $o$ whose divisors in  $\widehat{S}_G $ display the affine $A$-$D$-$E$ diagram property, we shall exclude the $A_1$ case, of which one can easily work out the solution.
\begin{theorem} \label{thm:uniF}
For $G \neq A_1$, let $F$ be the function defined in ${\rm Theorem \ \ref{thm:BrI} }$. Then $F$ is the unique local function of $S_G$ (up to the multiplication of local units) near the singular point $o$ with the affine $A$-$D$-$E$ diagram property, i.e., if $f$ is a local analytic function of $S_G$ near $o$ so that its divisor ${\rm div}( \pi^* f)$ in $\widehat{S}_G$ is described by the affine diagram corresponding to $G$, then  $f = u F$ for some unit $u$ as local functions of the analytic germ $( S_G , $o$ )$. (For the $D_r$ cases, the scalar $c$ in the definition of $F$ is uniquely determined by the given $f$). 
\end{theorem}
In this article, the study will be discussed in the domain of orbifold geometry. It is known that the groups $G$ of type $D$ or $E$ all have minus-identity in the center subgroup. By first blowing up $\CZ^2$ at the origin, the construction of the minimal resolution of $S_G$ can be reduced to that of $A$-type groups, where explicit coordinate systems are available by techniques in toric geometry. Using the explicit forms of invariants $X, Y, Z$, one can carry out the necessary calculations related to the divisors involved in the above two theorems, then make the connection with the affine $A$-$D$-$E$ diagram to achieve the results.

The paper is organized
as follows.  In Section 2, we briefly discuss the explicit forms of $A$-$D$-$E$ Kleinian groups and the expressions of invariants $X, Y, Z$ in Table 1, which our later calculations in the proof of Theorem \ref{thm:BrI} and \ref{thm:uniF} will base on. The verification of Theorem \ref{thm:BrI} will be given in the next two sections. Section 3 deals with the $A$-$D$ type groups $G$, where the explicit coordinates in toric geometry for the minimal resolution of related $A$-type groups are used, and some extract consideration of a certain $\ZZ_2$-symmetry is needed in the derivation. In Section 4, we consider the $E$-type groups $G$ in Theorem \ref{thm:BrI} by using the orbifold of $S_{A_1}$ quotiented by $\overline{G} \ (=G/A_1)$ to study the structure of $\widehat{S}_G$, then the conclusion will follow from calculations of $A$-type group action around three elements in the exceptional curve of $S_{A_1}$. In Section 5, we provide the proof of Theorem \ref{thm:BrI}. Based on the explicit forms of local $G$-invariant functions near the singular point $o$, all local functions in $S_G$ with the affine $A$-$D$-$E$ diagram property in Theorem \ref{thm:uniF} can be related to the function $F$ in Theorem \ref{thm:BrI} in a precise manner through various refine models associated to $\widehat{S}_G$. Finally we give some concluding remarks in Section 6.    

{\bf Convention}. 
To present 
our work, we prepare some notations. In this
paper, 
$\ZZ, \CZ$ will denote 
the ring of integers, complex numbers
respectively, $\PZ^1$= the complex projective line, and ${\rm i} = \sqrt{-1}$. For a positive integer $N$, we denote by $\ZZ_N$ the quotient ring $\ZZ/N\ZZ$. 

\section{Invariants of Kleinian Singularities} 
For the need of computations later used in this paper, we shall fix a $\CZ^2$-representation for each  group $G$ in Table 1, and denote the coordinates of $\CZ^2$ by $(Z_1, Z_2)$. For simplicity, the group $G$ will be denoted again by the corresponding symbol, $A,D$, $E$ if no confusion could arise.

Denote $\omega_N := e^{\frac{2 \pi {\rm i}}{N}}$, and define the following matrices in ${\rm SU}_2$:
$$
\sigma_N = \left( \begin{array}{cc} \omega_N  & 0 \\
                                    0 & \omega_N^{-1}  \end{array} \right) \ \ \ \ {\rm for} \ \  N \geq 2 \  \ , 
$$
and
$$
\tau = \left( \begin{array}{cc} 0 & 1 \\
                                   -1 & 0   \end{array} \right) \ , \ \  \mu = \frac{1}{\sqrt{2}}\left( \begin{array}{cc} \omega_8^7 , & \omega_8^7   \\
                                   \omega_8^5 ,  & \omega_8^7   \end{array} \right) \ , \ \ \kappa = \frac{1}{\sqrt{5}}\left( \begin{array}{cc} \omega_5^4 - \omega_5 ,
& \omega_5^2  - \omega_5^3   \\
                                 \omega_5^2  - \omega_5^3 ,  & \omega_5  - \omega_5^4 
 \end{array} \right)
$$
We have the relations: $\tau \sigma_N  = \sigma_N^{-1} \tau$, $\tau^2= \mu^3 = \kappa^2 = -1$, $\mu \tau = - \sigma_4 \tau \mu $ and $\mu \sigma_4 = \tau^{-1} \mu$.
The finite subgroup $G$ of ${\rm SL}_2(\CZ)$ can be represented by the following forms:
\be
A_r = \langle \sigma_{r+1} \rangle \ \ , \ \ D_r  = \langle \sigma_{2r-4} , \tau \rangle \ \ , \ \
E_6  = \langle \sigma_4 , \tau , \mu \rangle \ \ , \ \
E_7 = \langle \sigma_8 , \tau , \mu \rangle \ \ , \ \
E_8  = \langle \sigma_{10} , \kappa \rangle \ \ ,
\ele(ADE)
(see, e.g., \cite{Sl} ). The invariants $X, Y, Z$ in Table 1 are given by 
\bea(llll)
A_r :& X = Z_1Z_2 , & Y = Z_1^{r+1} , & Z = Z_2^{r+1} \ ; \\
D_r :& X = Z_1^2Z_2^2 ,& Y = \frac{1}{2}(Z_1^{2r-4} + Z_2^{2r-4}) , & Z = \frac{1}{2}Z_1Z_2(Z_1^{2r-4} - Z_2^{2r-4}) \ ; \\[1mm]
E_6 :& X= f_6(Z_1, Z_2), & Y = \frac{-1}{3\sqrt[3]{4}}f_8(Z_1, Z_2), & Z= \frac{1}{6\sqrt{3}} f_{12}(Z_1, Z_2) ; \\
E_7 :& X= \frac{-1}{\sqrt[3]{3}}f_8(Z_1, Z_2), & Y = -6 f_6 (Z_1, Z_2)^2, & Z=  {\rm i} \sqrt{2}  f_6(Z_1, Z_2) f_{12}(Z_1, Z_2) ;  \\
E_8 :& X= -\sqrt[5]{1728}F_{12}(Z_1, Z_2), & Y = \frac{1}{121} F_{20}(Z_1, Z_2) , & Z=  \frac{1}{20} {\rm Det(Jacobi}(F_{12}, F_{20})),
\elea(invt)
where $f_j, F_k$'s are polynomials in $Z_1, Z_2$ (see \cite{Kl} or \cite{H}) given by 
$$
\begin{array}{ll}
f_6 (Z_1, Z_2) = Z_1Z_2 (Z_1^4-Z_2^4), & f_8(Z_1, Z_2) = Z_1^8 + 14 Z_1^4Z_2^4 + Z_2^8 , 
\\
f_{12}(Z_1, Z_2)= Z_1^{12}-33Z_1^8Z_2^4-33Z_1^4Z_2^8+Z_2^{12}, &  \\
F_{12}(Z_1, Z_2)= Z_1 Z_2(Z_1^{10}+11 Z_1^5 Z_2^5 -Z_2^{10}), & F_{20}(Z_1, Z_2) ={\rm Det(Hess}(F_{12})) .
\end{array}
$$
For a group $G$ in (\req(ADE)), it induces a subgroup $\overline{G}$ in ${\rm PSL_2}(\CZ)$ with $|\overline{G}| = |G|/2$ except $G= A_{2k}$, in which case $|\overline{G}|=|G|$. The projective automorphism group $\overline{G}$ acts on $\PZ^1$ with the ratios $[Z_1 , Z_2]$ as the homogenous coordinates.  The $\overline{G}$-quotient of $\PZ^1$ gives rise to a morphism of $\PZ^1$,
\be
\varphi : \PZ^1 \longrightarrow \PZ^1 \ ( = \PZ^1/\overline{G} ) \ , \ \ [Z_1 , Z_2] \mapsto  [\varphi_1(Z_1, Z_2) , \varphi_2 (Z_1, Z_2) ] \ ,
\ele(GP1)
where $\varphi_1, \varphi_2$ are two homogenous polynomials of degree $|\overline{G}|$. For $E$-type groups, it is known that $\varphi_1, \varphi_2$ can be represented by the following invariants, and the map (\req(GP1)) has three critical values, $v_1, v_2, v_\infty$. All critical points over a critical value $v_j$ have the same branched index, denoted by $b_j$, equal to the order of their  isotropic subgroups in $\overline{G}$: 
\bea(|l | c | c | c |)
\hline 
& [\varphi_1 , \varphi_2] & {\rm Critical \ Value :} \ v_1, v_2, v_\infty & {\rm Branched \ index : } \ b_1 , b_2 , b_\infty \\ \hline
E_6 & [ Z , X^2 ] & [{\rm i}  , 1] , \ [-{\rm i} , 1] , \ [1 , 0] & 3, \ 3, \ 2 \\[1mm]
E_7 & [ Y^2 , X^3 ] & [0 , 1], \ [1 , 1] , \ [ 1 , 0]& 4 , \ 2 , \ 3 \\[1mm]
E_8 &  [ Y^3 , X^5 ] & [0 , 1] , \ [-1  , 1 ] \ , [1 , 0 ] & 3 , \ 2  , \ 5 \\ \hline
\elea(EP1)

\section{The Affine Diagram of Type A or D}
In this section, we are going to show Theorem \ref{thm:BrI} for groups $G$ of type $A$ or $D$. 

For $G= A_r$, one can use techniques in toric varieties \cite{DMO} or the
continued fraction method \cite{HH} to obtain 
the minimal resolution
$\widehat{S}_G$. Indeed, $\widehat{S}_G$ has an open cover consisting of open affine charts, $ {\cal U}_k
\ ( 0 \leq k \leq r )$, where ${\cal U}_k$ is isomorphic to $\CZ^2$ with the affine  
coordinates $(u_k, v_k)$ expressed by the monomials of $Z_i$'s : 
\be
{\cal U}_k \simeq \CZ^2 \ni (u_k, v_k) =
(Z_1^{k+1}Z_2^{-r+k} , Z_1^{-k}Z_2^{r+1-k} ) \ \ \ \ {\rm for } \ \ \ 0 \leq k \leq r \ .
\ele(uvk)
Denote by $\hat{o}_k$ the element in $\widehat{S}_G$
with the coordinate $u_k=v_k=0$. The exceptional
divisor in $\widehat{S}_G$
is given by $\ell_1+ \cdots + \ell_r$, where $\ell_j$ is  
the rational $(-2)$-curve joining
$\hat{o}_{j-1}$ and $\hat{o}_j$ with the  defining equations: $v_{j-1} = u_j = 0$. Inside $S_G$ and $\widehat{S}_G$, the algebraic torus $\CZ^{* 2}/G$ is considered as the same Zariski-open set via the morphism (\req(resl)). As a toric variety, the complement of $\CZ^{* 2}/G$ in $\widehat{S}_G$ consists of $r+2$ irreducible toric divisors, $c_0+ \sum_{j=1}^r \ell_j + c_{r+1}$. Here $c_0, c_{r+1}$ are the affine curve $\CZ$ defined by $u_0=0$ and $v_r=0$  respectively, corresponding to the proper transform of $\{Z_i=0 \}/G$ for $i=1, 2$. The structure of $\widehat{S}_G$ can be depicted in the following picture:
$$
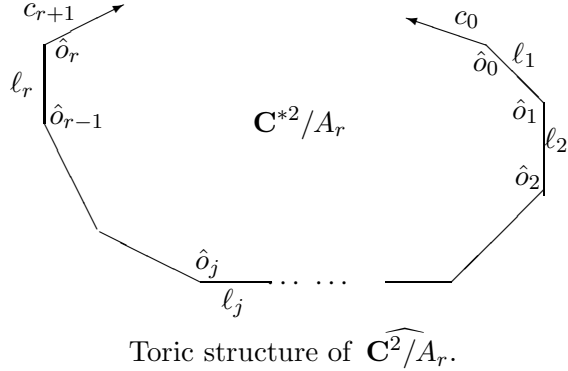
\put(73, -2){\shortstack{$c_0$ }}
\put(85, -10){\vector(-3, 1){30}}
\put(85, -10){\line(1, -1){22}}
\put(107, -32){\line(0, -1){35}}
\put(107, -65){\line(-1, -1){35}}
\put(72, -100){\line(-1, 0){25}}
\put(2, -103){\shortstack{$\cdots \ \cdots$ }}
\put(4, -100){\line(-1, 0){27}}
\put(-23, -100){\line(-2, 1){38}}
\put(-62, -80){\line(-1, 2){20}}
\put(-82, -40){\line(0, 1){30}}
\put(-82, -10){\vector(2, 1){30}}
\put(80, -20){\shortstack{$\hat{o}_0$ }}
\put(95, -38){\shortstack{$\hat{o}_1$ }}
\put(96, -63){\shortstack{$\hat{o}_2$ }}
\put(-80, -40){\shortstack{$\hat{o}_{r-1}$ }}
\put(-78, -15){\shortstack{$\hat{o}_r$ }}
\put(-25, -95){\shortstack{$\hat{o}_j$ }}
\put(95, -17){\shortstack{$\ell_1$ }}
\put(107, -49){\shortstack{$\ell_2$ }}
\put(-15, -110){\shortstack{$\ell_j$ }}
\put(-95, -28){\shortstack{$\ell_r$ }}
\put(-90, 1){\shortstack{$c_{r+1}$ }}
\put(-5, -43){\shortstack{$\CZ^{* 2}/A_r$ }}
\put(-50, -130){\shortstack{
Toric structure of $\widehat{\CZ^2/A_r}$. }}
$$ 
By (\req(uvk)), one has 
\be
Z_1 Z_2 = u_k v_k  \ , \ \ Z_1^{r+1} = u_k^{r-k+1}v_k^{r-k} \ , \ \ Z_2^{r+1} = u_k^k v_k^{k+1} \  \ \ \ {\rm for} \ \ \ 0 \leq k \leq r \ \ .
\ele(Arel)
Then using the expressions in (\req(invt)) for $A_r$, we have the following divisor-description in $\widehat{S}_{A_r}$:
\bea(ll)
{\rm div}(\pi^*X) = c_0+ \sum_{j=1}^r \ell_j + c_{r+1} ,& \\
 {\rm div}(\pi^*Y) = (r+1)c_0 + \sum_{j=1}^r (r+1-j) \ell_j , &
{\rm div}(\pi^*Z) =  \sum_{j=1}^r j \ell_j +(r+1)c_{r+1} \ .
\elea(Adiv)
Then $\pi^*X$ gives rise to a regular function of $\widehat{S}_{A_r}$, whose divisor is represented by the diagram $\widetilde{A}_r$ with $\oplus$ corresponding $c_0+c_{r+1}$, which is the proper transform of $(X=0)$ in $S_{A_r}$. Hence we obtain Theorem \ref{thm:BrI} for $G= A_r$.

For $G= D_r$ ($r \geq 4$), we denote $h = 2r-5$, and $H =$ the group $A_h$.  By (\req(ADE)), $H $ is an index $2$ normal subgroup of $G$, and  one has the naturally induced morphism of orbifolds: $S_H \longrightarrow S_G$. The element $\tau$ in $G$ defines the linear automorphism of $\CZ^2$ sending $(Z_1, Z_2)$ to $(Z_1^*, Z_2^*) := (-Z_2, Z_1)$, hence it induces an order 2 automorphism of $S_H$, by which one has the volume-form preserving automorphism $\tau^\dagger$ of $\widehat{S}_H$. The minimal resolution $\hat{S}_G $ of $S_G$ can be identified with the minimal resolution of $\widehat{S}_H/\tau^\dagger$. For the affine charts ${\cal U}_k \ (0 \leq k \leq h)$ of $\widehat{S}_H$ in (\req(uvk)), $\tau^\dagger$ interchanges ${\cal U}_k$ and ${\cal U}_{h-k}$ with the following coordinate-expression:
$$
{\cal U}_k \ni (u_k, v_k) \mapsto (u_k^*, v^*_k) = (-1)^k(- v_{h-k},  u_{h-k}) \ , \ \ {\rm where} \ (u_{h-k}, v_{h-k}) \in {\cal U}_{h-k} \ .
$$
Therefore, the fixed-point set of $\tau^\dagger$ in $\widehat{S}_H$ consists of two elements, denoted by $\alpha_1, \alpha_2$. Indeed by $k=r-3$ in the above relations, $\tau^\dagger$ is described by
$$
(u_{r-3}, v_{r-3}) \mapsto 
(-1)^r ( v_{r-2}, - u_{r-2}) = (-1)^r ( u_{r-3}^{-1}, - u_{r-3}^2 v_{r-3}) \ .
$$
The elements $\alpha_1, \alpha_2$ are defined by the equations: $v_{r-3}=0, u_{r-3}^2 = (-1)^r$. This implies that both $\alpha_i \in \ell_{r-2}$, and  $\tau^\dagger$ acts locally as $A_1$ near $\alpha_i$ in $\widehat{S}_H$. Let $\widetilde{S}_H$ be the blow-up of $\widehat{S}_H$ centered at $\{ \alpha_1, \alpha_2 \}$, and $\widetilde{e}_1, \widetilde{e}_2$ the corresponding exceptional (-1)-curves in $\widetilde{S}_H$. The automorphism $\tau^\dagger$ of $\widehat{S}_H$ can be lifted to an involution $\widetilde{\tau}$ of $\widetilde{S}_H$ with $\widetilde{e}_1 \cup \widetilde{e}_2$ as the fixed locus. Then the quotient space $\widetilde{S}_H/\widetilde{\tau}$ is a smooth surface which can be identified with the resolution space $\widehat{S}_G$ of $S_G$,
\be
\wp : \ \widetilde{S}_H \ \longrightarrow \widehat{S}_G = \widetilde{S}_H/\widetilde{\tau} \ .
\ele(Dcov)
The exceptional divisors in $\widehat{S}_G$ are the $\wp$-image of $\widetilde{e}_i$'s and $ \ell_j$'s, described as follows:
$$
e_i := \wp (\widetilde{e}_i) \ \ (i=1, 2) \ , \ \ d_j := \wp (\ell_j ) = \wp (\ell_{2r-4-j}) \ \ (1 \leq j \leq r-3) \ , \ \ d_{r-2}:= \wp (\ell_{r-2} ) \ .
$$
Note that for a function $f$ of $S_G$, one can determine the structure of ${\rm div}(\pi^*f)$ in $\widehat{S}_G$ by examining that of ${\rm div}((\pi \wp)^*f)$ in $\widetilde{S}_H$ through the (degree 2) morphism $\wp$ in (\req(Dcov)). Since $\tau^\dagger$ interchanges ${\cal U}_k$ and ${\cal U}_{h-k}$, the expressions of the function $f$ on the first half of open charts in $\widehat{S}_H$, ${\cal U}_j $ for $ 0 \leq j \leq r-3$, will explicitly display the structure of ${\rm div}(\pi^*f)$ in $\widehat{S}_G$. 
By (\req(Arel)) for $A_h$, one can express the following  functions of $S_G$ in the coordinates $(u_j, v_j)$ of ${\cal U}_j $ for $ 0 \leq j \leq r-3$,
$$
2Y \ ( = Z_1^{h+1} + Z_2^{h+1}) = u_j^j v_j^{j+1}( u_j^{2r-4-2j}v_j^{2r-6-2j} + 1)
\ , \ \
X \ (= Z_1^2 Z_2^2 ) = u_j^2 v_j^2 . 
$$
Hence for $c \in \CZ \setminus \{ 0 \}$, one has 
\be
X + c Y = \left\{ \begin{array}{ll} 
v_0(u_0^2 v_0 + \frac{c}{2}  u_0^{2r-4}v_0^{2r-6} + \frac{c}{2} ) & {\rm for} \ j = 0 ; \\
u_1 v_1^2(u_1 + \frac{c}{2} u_1^{2r-6}v_1^{2r-8} + \frac{c}{2}) & {\rm for} \ j = 1 ; \\
u_j^2 v_j^2(1 + \frac{c}{2} u_j^{2r-j-6}v_j^{2r-j-7} + \frac{c}{2} u_j^{j-2} v_j^{j-1}) & {\rm for} \ 2 \leq j \leq r-3 .  \end{array} \right.
\ele(DF)

In the case $r=4$, (hence $h=3$), the elements $\alpha_1, \alpha_2$ in $\widehat{S}_H$ are given by $(u_1 , v_1)=(\pm 1 , 0)$.
By (\req(DF)) one has 
$$
X + c Y = v_0(u_0^2 v_0 + \frac{c}{2}  u_0^4 v_0^2 + \frac{c}{2} ) = \frac{c}{2} u_1 v_1^2( u_1^2 +\frac{2}{c}u_1 + 1) \ .
$$
When $c \neq \pm 1$, the zeros of the function $u_1^2 +\frac{2}{c}u_1 + 1$ consist of two disjoint $\CZ$-curves, denoted by $\widetilde{\varrho}_1 , \widetilde{\varrho}_2$. Note that $\widetilde{\varrho}_1 , \widetilde{\varrho}_2$ are disjoint with the elements $\alpha_1, \alpha_2$, and are permuted under $\widetilde{\tau}$. Therefore for $F = X + c Y \ ( c \in \CZ \setminus \{0, \pm 1 \})$, the divisor of $F$ in $\widetilde{S}_H$ near $(\pi \wp)^{-1}(0)$ is given by  ${\rm div}((\pi \wp)^*F) = \ell_1 + \ell_3 + 2 \ell_2 + 2(\widetilde{e}_1 + \widetilde{e}_2) + \widetilde{\varrho}_1 + \widetilde{\varrho}_2$. Hence in $\widehat{S}_G$ near the exceptional set, one obtains 
\be
{\rm div}(\pi^*F) = 2 d_2 + d_1 + e_1 + e_2 + \varrho  \ ,
\ele(D4) 
with the $\widetilde{D}_4$-configuration,
where $\varrho = \wp (\widetilde{\varrho}_1 )= \wp (\widetilde{\varrho}_2)$, the $\CZ$-curve corresponding to $\oplus$ of the diagram.

For $r \geq 5$, $F = X + c Y$ with $c \neq 0$. Denote by $\widetilde{\varrho}$ the curve in $\widehat{S}_H$ defined by  $u_1 + \frac{c}{2} u_1^{2r-6}v_1^{2r-8} + \frac{c}{2} =0$  in  the coordinates $(u_1, v_1)$ of ${\cal U}_1$. As both $\widetilde{\varrho}$ and $\tau^\dagger(\widetilde{\varrho})$ are disjoint with $\ell_{r-2}$, and $\alpha_1, \alpha_2 \in \ell_{r-2}$, $\widetilde{\varrho}$ and $\tau^\dagger(\widetilde{\varrho})$ can also be regarded as two affine curves in $\widetilde{S}_H$, denoted by  $\widetilde{\varrho}_1$, $\widetilde{\varrho}_2$ respectively. By (\req(DF)), the divisor  of the function $F$ in $\widetilde{S}_H$ near $(\pi \wp)^{-1}(0)$ is expressed by ${\rm div}((\pi \wp)^*F) = \ell_1 + \ell_h + 2 \sum_{j=2}^{h-1} \ell_j + 2(\widetilde{e}_1 + \widetilde{e}_2) + \widetilde{\varrho}_1 + \widetilde{\varrho}_2$. Therefore in $\widehat{S}_G$ near the exceptional set, one obtains the $\widetilde{D}_r$-configuration:
\be
{\rm div}(\pi^*F) =  d_1 + 2 \sum_{j=2}^{r-2} d_j +  e_1 + e_2 + \varrho  \ ,
\ele(Dr)
where $\varrho = \wp (\widetilde{\varrho}_1 )= \wp (\widetilde{\varrho}_2)$, the open affine curve corresponding to $\oplus$ of the diagram. It is easy to see that the affine curve $\varrho$ in $\widehat{S}_{D_r}$ is the proper transfer of $(F=0)$ in $S_{D_r}$ for $r \geq 4$. This completes the proof of Theorem \ref{thm:BrI} for $D$-type groups $G$.
\par \noindent
{\bf Remark.} In the above proof of Theorem \ref{thm:BrI} for $G$=$D_r$, the open curve $\varrho$, corresponding to $\oplus$ of the affine $D_r$-diagram, intersects (normally) only $d_2$ among all the exceptional (-2)-$\PZ^1$ divisors of $\widehat{S}_G$. Note that the invariant function $F$ depends on a parameter $c$ in this case. The choice of $c$ corresponds to  the intersecting point $\varrho \cdot d_2$ in $d_2 \setminus \bigcup \{  e_1 , e_2 , d_j \ (j \geq 2 ) \}$ in the one-one correspondence manner.
$\Box$ 

\section{The Affine Diagram of Type E}
In this section, we will show Theorem \ref{thm:BrI} for $G= E_6, E_7, E_8$, where the function $F$ is defined by $F = X$. 

Note that $A_1$ is the central subgroup of $G$. The blow-up of $\CZ^2$ at the origin, $\widehat{\CZ^2}$, has the (-1)-$\PZ^1$ exceptional curve with the ratios $[Z_1, Z_2]$ as its homogenous coordinates. The $A_1$-action on $\CZ^2$ induces a $\ZZ_2$-action of $\widehat{\CZ^2}$ with the exceptional curve  as its branched locus. The $\ZZ_2$-quotient of $\widehat{\CZ^2}$ is a smooth surface which can be identified with the minimal resolution of $S_{A_1}$, denoted by $M \ (:= \widehat{S}_{A_1})$ in this section. We denote the (-2)-exceptional curve in $M$ by $\widetilde{e}$, which is isomorphic to the exceptional curve in $\widehat{\CZ^2}$, hence one may again use $[Z_1, Z_2]$ as the  homogenous coordinates of $\widetilde{e}$. The action of $\overline{G} \ (= G/A_1)$ on $S_{A_1}$ induces a $\overline{G}$-action on $M$ preserving the global holomorphic volume-form. Then we have the following morphism induced from $S_{A_1}$ to its $\overline{G}$-quotient, 
\be
\Phi : M = \widehat{S}_{A_1} \ \longrightarrow S_G \ (= S_{A_1}/\overline{G} ) \ . 
\ele(Phi)
The minimal resolution $\widehat{M/\overline{G}}$ of $M/\overline{G}$ is isomorphic to $\widehat{S}_G$, and we shall make the identification, $
\widehat{S}_G = \widehat{M/\overline{G}}$, for the rest of this section.  
We are going to determine the structure of ${\rm div} (\pi^* X)$ in $\widehat{S}_G$ by examining the divisor ${\rm div}(\Phi^* X)$ in $M$ under the $\overline{G}$-action. 

The exceptional curve $\widetilde{e}$ of $M$ is $\overline{G}$-stable, and it contains all the $\overline{G}$-fixed points of $M$: $M^{\overline{G}} = \widetilde{e}^{\overline{G}}$. With $[Z_1, Z_2]$ as the coordinates of $\widetilde{e}$, the quotient map from $\widetilde{e}$ to $\widetilde{e}/\overline{G}$ can be represented by the morphism $\varphi$ in (\req(GP1)) with the property (\req(EP1)), and we shall make this identification hereafter for simplicity of notations. Hence  
$M^{\overline{G}}$ consists of all critical points of the morphism $\varphi$, $M^{\overline{G}} = \bigcup_{j=1, 2, \infty} \varphi^{-1}(v_j)$. Furthermore for each  
$p \in \varphi^{-1}(v_j)$, the isotropic subgroup $\overline{G}_p$ at $p$ is an order $b_j$ cyclic subgroup of $\overline{G}$ with the $A_{b_j-1}$-action near $p$ in $M$, where $b_j$ is the branched index in (\req(EP1)). By (\req(invt)), the linear factors (in the variables $Z_1, Z_2$) of the invariant polynomial $X$ for each $E$-type group are all distinct with the multiplicity one. The $[Z_1, Z_2]$-ratios of these linear factors are indeed all the elements of $\varphi^{-1} ( v_{\infty} )$. For each $q \in \varphi^{-1} ( v_{\infty} )$, its corresponding linear form appeared in the factorization of $X$ gives rise to an irreducible divisor in $M$, denoted by $\widetilde{d}_q$, which intersects $\widetilde{e}$ transversely. It is not hard to see that the divisor in $M$ induced by the function $X$ is given by
\be
{\rm div} (\Phi^* X) = m \widetilde{e} + \sum_{q \in \varphi^{-1} ( v_{\infty} )} \widetilde{d}_q \ \ , \ \ \ \ \ \ \ {\rm where} \ m:= \frac{{\rm deg} (X)}{2} \ .
\ele(PhiX) 
For a $\overline{G}$-fixed point $p \in \varphi^{-1}(v_j)$, there exists a local analytic coordinates $(z_1, z_2)$ for a $\overline{G}_p$-invariant neighborhood $U_p$ centered at $p$ in $M$, such that $\widetilde{e}$ is locally defined by $z_1 = 0$ and the local $A_{b_j-1}$-action of $\overline{G}_p$ is generated by $(z_1, z_2) \mapsto (\omega_{b_j} z_1, \omega_{b_j}^{-1} z_2  )$; furthermore for $j = \infty$, the divisor $\widetilde{d}_p$ is locally defined by $z_2 = 0$. In terms of $(z_1, z_2)$, the function $\Phi^* X$ near $p \in \varphi^{-1}(v_j)$ has the following local description:
\be
\Phi^* X = \left\{ \begin{array}{ll} z_1^m & {\rm for} \ j = 1 , 2 ; \\
z_1^m z_2   & {\rm for} \ j = \infty \ .
\end{array} \right.
\ele(Eleq)
Note that $m = b_\infty +1$, and $m$ is divisible by $b_1$ and $b_2$ by (\req(EP1)). 

We now consider the orbifold $M/\overline{G}$. The singular set of  $M/\overline{G}$ consists of three elements, corresponding to the critical values of $\overline{G}$-quotient map of $\widetilde{e} \ ( = \PZ^1 )$ in (\req(GP1)), denoted again by $v_1, v_2, v_\infty$. The local orbifold structure of $M/\overline{G}$ near $v_j$ is isomorphic to $U_p/\overline{G}_p$ for any $p \in \varphi^{-1}(v_j)$. Using the local coordinates $(z_1, z_2) \in \CZ^2$ of $U_p$, one has the isomorphism of germs of analytic spaces:
$$
(M/\overline{G} , v_j ) \ \simeq \ (U_p/\overline{G}_p , [p] ) \ = ( \CZ^2/A_{b_j-1} , [\vec{0}] ) \ ,
$$ 
hence the structure of $\widehat{M/\overline{G}}$ over $M/\overline{G}$ is described by $\widehat{S}_{A_{b_j-1}}$ ($j=1, 2, \infty$) near the exceptional set. In order to determine the divisor of the function $X$ in $\widehat{M/\overline{G}}$, one needs only to consider the divisor-expression of $X$ in the minimal resolution $\widehat{U_p/\overline{G}_p}$ ( $\simeq \widehat{S}_{A_{b_j-1}}$ locally) of $U_p/\overline{G}_p$, for $p \in \varphi^{-1}(v_j)$. With $(Z_1, Z_2) = (z_1, z_2)$ in (\req(Arel)) for $r=b_j$, one has the divisor-description of the following local functions:
$$
\begin{array}{lll}
{\rm div}(z_1^{b_j}) &= b_j e (p) + \sum_{k=1}^{b_j-1} (b_j-k) \ell_{p, k}  & {\rm for \ all} \ j \ , \\
{\rm div}(z_1z_2) &= e (p) +d (p)   + \sum_{k=1}^{b_\infty-1} \ell_{p, k} & {\rm for} \ j = \infty \ , 
\end{array}
$$
where $\ell_{p, k}$'s are the (-2)-$\PZ^1$ exceptional curves, and $e(p), d(p)$ are the open curves in $\widehat{U_p/\overline{G}_p}$ corresponding to the proper transfer of $\widetilde{e}/\overline{G}_p$, $\widetilde{d}_p/\overline{G}_p$ respectively. By (\req(Eleq)), the function in $\widehat{U_p/\overline{G}_p}$ induced by $X$ has the following divisor-expression:
$$
{\rm div}(X) = \left\{ \begin{array}{ll} 
m e (p) + \sum_{k=1}^{b_j-1} (m-\frac{mk}{b_j}) \ell_{p, k}  &{\rm for} \ j = 1 , 2 , \\
m e (p) + d(p) + \sum_{k=1}^{b_\infty-1} (m-k) \ell_{p, k}  &{\rm for} \ j = \infty \ .  \end{array} \right.
$$
Note that the function of $M$ induced by $X$, i.e. $\Phi^* X$, is $\overline{G}$-invariant, so is  ${\rm div}(\Phi^* X)$. Furthermore, the local resolution spaces $\widehat{U_p/\overline{G}_p}$ for $p \in \varphi^{-1}(v_j)$ are all identified in $\widehat{M/\overline{G}}$ for each $j$. Therefore for $1 \leq k \leq b_j-1$,  all the curves $\ell_{p, k}$  ( $ p \in \varphi^{-1}(v_j)$) are the same (-2)-$\PZ^1$ exceptional curve in $\widehat{M/\overline{G}}$, denoted by $\ell_{v_j, k}$. In $\widehat{M/\overline{G}}$, we consider the curves, $e =$ the proper transfer of $\widetilde{e}/\overline{G} \ ( \subset M/\overline{G})$, and $d =$ the proper transfer of ${\rm div}(X) \ ( \subset S_G )$. Then $e$ is a (-2)-$\PZ^1$ curve and $d$ is an affine curve in $\widehat{M/\overline{G}}$. Through the $\overline{G}$-structure of $M$, the curves $e(p), d(p)$ in $\widehat{U_p/\overline{G}_p}$ are open subsets of $e, d$ respectively. By the local results on the structure of $\widehat{U_p/\overline{G}_p}$, the exceptional divisors of $\widehat{M/\overline{G}} \ (= \widehat{S}_G)$ over $S_G$ is given by 
$ e + \sum_{j=1, 2, \infty} \sum_{k=1}^{b_j-1} \ell_{v_j, k}$ with the corresponding $E$-type configuration, and the divisor of $\pi^* F$ in $\widehat{M/\overline{G}}$ for $F=X$ is expressed by
\be
{\rm div}(\pi^* F) = m e + d + \sum_{k=1}^{b_\infty -1} (m-k) \ell_{v_\infty, k} + \sum_{j=1}^2 \sum_{k=1}^{b_j-1} (m-\frac{mk}{b_j}) \ell_{v_j, k} \ .
\ele(EFd) 
By which, one obtains the affine $E$-diagram property of ${\rm div}(\pi^* X)$ for each case, where $d$ corresponds to $\oplus$ of the diagram, which is the proper transform of $(X=0)$ in $S_G$. This completes the proof of Theorem \ref{thm:BrI} for groups $G$ of type $E$.

\section{The Uniqueness of Local Invariants with the Affine A-D-E Diagram Property}
Now we are going to provide a proof of Theorem \ref{thm:uniF}. For convenience of the discussion, throughout this section a local function of $S_G$ always means a function near the singular point $o$, and the divisor of a function in some variety over $S_G$ always denote  its divisor near the exceptional set over $o$. The function  $F$ defined in Theorem \ref{thm:BrI} is now considered as a local function of $S_G$, and $f$ will be a local function with the affine $A$-$D$-$E$ diagram property in Theorem \ref{thm:uniF}.

For $G= A_r \ (r \geq 2)$, the divisor of $F$ in $\widehat{S}_G$ is given by (\req(Adiv)) for $F=X$,  
$$
{\rm div}(\pi^*F) = c_0+ \sum_{j=1}^r \ell_j + c_{r+1} \ .
$$
By the $A_r$-diagram property of $f$, the divisor of $f$ in $\widehat{S}_G$ is expressed by
\be
{\rm div}(\pi^*f) = c_0^\prime+ \sum_{j=1}^r \ell_j + c_{r+1}^\prime \ ,
\ele(Afd)
where $c_0^\prime, c_{r+1}^\prime$ are two disjoint open affine curves in $\widehat{S}_G$ near the exceptional set, which interest only $\ell_1, \ell_r$ respectively among all $\ell_j$'s. Note that $\ell_1 \neq \ell_r$ by $r \geq 2$, and no relation exists between two normal-crossing-point sets, $\{ c_0^\prime \cdot \ell_1,  c_{r+1}^\prime \cdot \ell_r \}$ and $\{ c_0 \cdot \ell_1,  c_{r+1} \cdot \ell_r \}$, at this moment. By the relation of invariants, $X^{r+1} = YZ$, one can express the local function $f$ in the form, $f = \beta_1 X^{n_1} + \beta_2 Y^{n_2} + \beta_3 Z^{n_3}$ for some positive integers $n_j$ and local functions $\beta_j$ in $S_G$ such that $\beta_k$ is either zero or a local unit for $k=2, 3$.
If both $\beta_2$ and $\beta_3$ are units, by (\req(Adiv)) the support of ${\rm div}(\pi^*f)$ is contained in the exceptional set $\bigcup_{j=1}^r \ell_j$, a contradiction to the description of ${\rm div}(\pi^*f)$ in (\req(Afd)). Replacing $\bigcup_{j=1}^r \ell_j$ by $c_0 \cup \bigcup_{j=1}^r \ell_j$ or $c_{r+1} \cup \bigcup_{j=1}^r \ell_j$ in the previous  argument, one can also rule out the cases when exact one of $\beta_2, \beta_3$ is a unit. Hence $\beta_2 = \beta_3 =0$, then by (\req(Adiv)), one has $n_1=1$ and $\beta_1$= a unit, The result for $G= A_r$ follows immediately.

For $G=D_r$, the divisor of $F$ in $\widehat{S}_G$ is given by (\req(D4)) or (\req(Dr)), where $\varrho$ corresponds to $\oplus$ of the affine $D_r$-diagram. Let $\varrho^\prime$ be the irreducible component in ${\rm  div}(\pi^* f)$ corresponding to $\oplus$ of the affine diagram. Then $\varrho^\prime$ intersects only one (-2)-$\PZ^1$ curve, $d_2$, in the exceptional set, and we denote the intersecting point $d_2 \cdot \varrho^\prime$ by $p$. By the remark in Section 3, there is a unique complex number $c$ as the scalar in the definition of $F$ such that $\varrho \cdot d_2$ is equal to the element $p$. As both ${\rm div}( \pi^* F)$ and ${\rm div}( \pi^* f)$ are represented by the same affine diagram configuration, and $\varrho, \varrho^\prime$ intersect $d_2$ normally at the same point $p$, the support of the divisor ${\rm div}( \pi^*( f / F ))$ is disjoint with the exceptional set $\pi^{-1} (o)$. This implies that $ f / F $ gives rise to a unit $u$ of $S_G$, hence follows the result.

For an $E$-type group $G$, we shall identify $\widehat{S}_G$ with the minimal resolution $\widehat{M/\overline{G}}$ of $M/\overline{G}$ via the map $\Phi$ in (\req(Phi)) as in Section 4. The divisor of $F$ in $\widehat{S}_G$ is given by (\req(EFd)), where $d$ denotes the curve corresponding to $\oplus$ of the affine diagram corresponding to $G$. And the divisor of $f$ in $\widehat{S}_G$ is expressed by 
$$
{\rm div}(\pi^* f) = m e + d^\prime + \sum_{k=1}^{b_\infty -1} (m-k) \ell_{v_\infty, k} + \sum_{j=1}^2 \sum_{k=1}^{b_j-1} (m-\frac{mk}{b_j}) \ell_{v_j, k} \ ,
$$
where $d^\prime$ corresponds to $\oplus$  of the affine diagram. Both $d$ and $d^\prime$ are open curves meeting (transversely) only with $\ell_{v_\infty, b_\infty -1}$ among all irreducible exceptional curves in $\widehat{S}_G$. Consider the divisor of $f$ in $M$ via the morphism (\req(Phi)). By the structure of ${\rm div}(\pi^* f)$ and ${\rm div}(\pi^* F)$,   the construction of the resolution space $\widehat{M/\overline{G}}$ enables one to conclude the following parallel version of (\req(PhiX)) for the divisor of $f$ in $M$,   
$$
{\rm div} (\Phi^* f) = m \widetilde{e} + \sum_{q \in \varphi^{-1} ( v_{\infty} )} \widetilde{d}^\prime_q \ \ , 
$$ 
where $\widetilde{d}^\prime_q$ is an open curve meeting $\widetilde{e}$ transversely at $q$ for all $q \in \varphi^{-1} ( v_{\infty} )$. Note that the blow-up of $\CZ^2$ ( with the coordinates $(Z_1, Z_2)$) is the double cover of $M$ with the branched locus $\widetilde{e}$, of which the ratios $[Z_1, Z_2]$ serve as the homogenous coordinates. The set $\varphi^{-1} ( v_{\infty} )$ consists of the ratios associated to all linear factors of the invariant polynomial $F \ (= X)$. Hence as local $G$-invariant functions of $\CZ^2$ near the origin, $f/ F$ is a local $G$-invariant unit. Then follows the result of Theorem \ref{thm:uniF}.

\section{Concluding Remarks}
In this work, we have provided a quantitatively geometrical interpretation of branching indices attached to an affine $A$-$D$-$E$ Coxeter-Dynkin diagram by the divisor theory in the minimal resolution $\widehat{S}_G$ of a Kleinian orbifold $S_G$. Our approach has been based on the orbifold aspect of Kleinian singularities by exploring the explicit form of $G$-invariant polynomials for the associated orbifold $S_G$, instead of working on its defining equation (which is a consequence of the invariants' form), an algebraic geometry method commonly practiced in the study of some general surface singularities.   
By examining the expressions of $G$-invariants associated to $S_G$, a particular function $F$ naturally arises in the context so that its divisor in $\widehat{S}_G$ reveals the full data of the corresponding affine diagram, including the branching indices and $\oplus$. This can be viewed as an affine quantitative version of elliptic fibration near certain types of simple singular fibers in Kodaira's elliptic surface theory. 
We have obtained two main results, Theorem \ref{thm:BrI} and \ref{thm:uniF}, stated in Section 1. The justification of Theorem \ref{thm:BrI} was explained in Section 3 and 4, by using the explicit form of invariant functions and techniques in toric geometry applying to resolution problems of orbifolds. Along the same path practiced in the previous sections, and appending some extract efforts to analyze the divisors of functions involved, we proved  Theorem \ref{thm:uniF} in Section 5 to conclude the uniqueness nature of the function $F$ naturally appeared in "geometrization" of affine $A$-$D$-$E$ diagrams. Note that the geometry corresponding to $\oplus$ of the affine diagram plays a significant role on the uniqueness property of $F$ in Theorem \ref{thm:uniF}. Of course, one may make a direct check without much effort on the affine $A$-$D$-$E$ diagram property of $F$ by working on the defining equation of $S_G$ through general techniques in surface theory. However, the symmetry nature of the orbifold $S_G$ in connection with group representations would be lost in the process, which has been one of our motivations for this work. Furthermore, in the study of Gorenstein orbifolds $\CZ^n/G$ for $n \geq 3$, the crepant resolution $\widehat{\CZ^n/G}$ of $\CZ^n/G$ , if it would exist, plays the role of the minimal resolution in Kleinian singularity theory. In particular when $n=3$, the existence of $\widehat{\CZ^3/G}$ has been known  for any finite subgroups $G$ in ${\rm SL}_3(\CZ)$ (see \cite{R} and references therein). Moreover the quantitative structure of $\widehat{\CZ^3/G}$ can be explicitly explored by methods in orbifold geometry, but the defining equation of $\CZ^3/G$ is indifferent for the construction of $\widehat{\CZ^3/G}$ in most cases. Even though  for finite groups $G$ in ${\rm SL}_3(\CZ)$ the McKay-correspondence-like dual-pairing between conjugacy classes and irreducible representations in accord to certain combinatoric data (now 2-dimensional simplicial complexes instead of graphs) has not yet been found, one could still pursuit the appropriate simplicial complex model through the structure of $\widehat{\CZ^3/G}$. By our initial success on the affine $A$-$D$-$E$ geometrical nature of Kleinian singularities, it is pertinent to ask the "affine"-version of the simplicial complex (with integers attached to vertices) through the crepant resolution of a Gorenstein orbifold, then search the relevant interpretation of group representations in the domain of orbifold geometry. A program along this line is now under consideration.

\section*{Acknowledgements}
The author thanks MSRI, Berkeley, California, where parts of this work were carried out in Spring 2004, for the kind hospitality. This paper is partially supported by NSC grant  92-2115-M-001-014 of Taiwan.

\end{document}